# Stability of Random Sums

S. Satheesh , N. Unnikrishnan Nair and E. Sandhya

*Department of Statistics, Cochin University of Science and Technology*
*Cochin-682 022. India.*
*Department of Statistics, Prajyoti Niketan College, Pudukkad*
*Trichur-680 301, **India.***

## Abstract

When the distribution of a random ($N$) sum of independent copies of a r.v $X$ is of the same type as that of $X$ we say that $X$ is N-sum stable. In this paper we consider a generalization of stability of geometric sums by studying distributions that are stable under summation w.r.t Harris law. We show that the notion of stability of random sums can be extended to include the case when $X$ is discrete. Finally we propose a method to identify the probability law of $N$ for which $X$ is N-sum stable.

**Keywords.** Probability generating function; Laplace transform; characteristic function; compounding; random sum; gamma; Linnik; Mittag-Leffler; stable; semi stable; semi-$\alpha$-Laplace; semi Mittag-Leffler; Bernoulli; Sibuya; geometric.

**AMS 2000 Subject Classification:** 60E05; 60E07; 60E10; 62E10; 62E15; 62E20.

## 1. Introduction

Let $X$, $X_1$, $X_2$, …. be non-degenerate independent and identically distributed (i.i.d) random variables (r.v) with a common characteristic function (CF) $\phi(u)$ and $N$ a positive integer valued r.v independent of $X$ with probability generating function (PGF) $P(t)$. Set $S_N = X_1 + X_2 +$ … $+ X_N$ which we refer to as the N-sum of $X_i$'s. When $cS_N$ and $X$ are identically distributed for some $c > 0$ ($S_N$ and $X$ are of the same type), or equivalently

$$P(\phi(cu)) = \phi(u), \quad \text{for all } u \in \mathbf{R}, \tag{1.1}$$

we say that the distribution of $X$ is stable under the operation of summation with respect to (w.r.t) the r.v $N$ ($X$ is N-sum stable). Notice that for N-sum stability we need relation (1.1) to be satisfied for some $c > 0$ only. Often $N$ is described as the compounding r.v.

Geometric summation schemes are the most studied and Kozubowski and Rachev (1999) have a comprehensive account of this (including their applications in different areas) and Kalashnikov (1997) discusses methods based on them. Some other works in this area are; Sandhya (1991a) characterizing semi-$\alpha$-Laplace laws by stability of geometric sums, Sandhya and Satheesh (1996) characterizing $\alpha$-Laplace (Linnik) laws by the class-L property among semi-$\alpha$-Laplace laws, Mohan (1992) discussing simultaneous stability of N-sum and random minimum to characterize exponential and geometric laws and Kozubowski and Podgorski (2000) discussing asymmetric Laplace laws. The geometric sum is related to Renyi's rarefaction (p-thinning) of renewal processes (Sandhya (1991b)) and has application in reliability theory also (Pillai and Sandhya (1996)). In more generality Sandhya (1991a, 1996), Gnedenko and Korolev (1996), and Kozubowski and Panorska (1996), have considered N-sum schemes to study N-infinitely divisible and N-stable laws.

But examples of non-geometric laws for *N* in the context appears to be rare and Sandhya (1991a, 1996) had discussed one for N-stable laws, motivated by Harris (1948).

When *X* has a discrete distribution we have the Poisson-sum characterization of modified logarithmic law resulting in negative binomial laws, in Milne and Yao (1989), and the known fact that geometric laws are closed under its own compounding. But no systematic discussion of N-sum stability of discrete laws is found in the literature. The reason is that we do not have a discrete analogue of the notion of distribution of the same type, which is at the heart of the concept of stability. Certain works that suggest the possibility of systematically defining N-sum stability of discrete laws as in (1.1), are Steutel and van Harn (1979), Devroye (1993), and Pillai and Jayakumar (1995).

Another problem is that of identifying the distribution of *N* that imparts stability of *X* (given) in N-sums. This is important as it could give some idea about the mechanism that generates the N-sum. In other words, the knowledge of the distribution of *X* characterizes that of *N*. This problem does not appear to have been discussed in the literature except as characterizations of geometric law in Cinlar and Agnew (1968) when *X* is exponential, Balakrishna and Nair (1997) when *X* is bivariate exponential and Sandhya (1996) when *X* is Mittag-Leffler.

With these three problems in mind this paper is organized in the following manner. In Section.2 we discuss stability of N-sums when *N* is Harris (a PGF discussed in Harris (1948)), which is a generalization of the geometric law. The N-sum stability of the generalized Linnik law of Devroye (1990), Erdogan and Ostrovskii (1998) is also presented here. The discussion of N-sum stability leads to simple proofs of results such as membership in class-L, unimodality, absolute continuity, etc. Developing the discrete analogue of distributions of the same type, we discuss stability of N-sums of discrete laws systematically in Section.3. Finally, in Section.4 we describe a method for identifying the PGF of *N* taking a clue from Sandhya (1991a, 1996).

We require the following class of continuous functions in our discussion.

$$\psi(u) = a\psi(bu) \quad \text{for all} \quad u \in \mathbf{R} \tag{1.2}$$

and some $a, b > 0$ with $\psi(o) = 0$ which has been discussed in the context of CFs of semi stable (SS) laws and regression equations (Kagan, et. al (1973), p. 9, 163, 323, 324), semi-$\alpha$-Laplace laws (Pillai (1985)) and integrated Cauchy functional equations (Pillai and Anil (1996)). It has been shown that $0 < b < 1 < a$ alone is possible and there exist a unique $\alpha > 0$ such that $ab^\alpha = 1$. Also when (1.2) is satisfied for two values of *b*, say $b_1$ and $b_2$ such that $\ln b_1/ \ln b_2$ is irrational then $\psi(u) = \lambda|u|^\alpha$, for some $\lambda > 0$. Further in $\exp\{-\psi(u)\}$ as the CF (of a SS law) $\alpha \in (0,2]$ and as the Laplace transform (LT) $\alpha \in (0,1]$.

## 2. Stability w.r.t Harris(*a,k*) Law

Let us begin with an aspect of N-sums implicit in equation (1.1).

**Lemma.2.1** In the set up of equation (1.1), $\phi(.)$ cannot have a real zero.

**Proof.** Notice that $P(s) = 0$ if and only if $s = 0$, since *P* is the PGF of a positive r.v. Again, we need consider only the case $c \in (0,1)$, because if $c > 1$ then a change of variable $v = cu$ will make (1.1) as $P(\phi(v)) = \phi(bv)$, where $b = 1/c < 1$. Now if $\phi(.)$ has a real zero say $\tau$, then $\phi(c\tau) = 0$ or $\phi(b\tau) = 0$ as the case may be because of the property of the PGF. Hence on iteration $\phi(c^n\tau) = 0$ for every $n > 1$ integer. Since $c \in (0,1)$ this will imply that

$$\underset{n\to\infty}{Lt}\ \phi(c^n\tau) = \phi(o) = 0,$$

also by the continuity of $\phi(s)$, which is a contradiction. Thus the assertion is proved.

**Remark.2.1** When $\phi(s)$ is a LT it is necessary that $c<1$. Because setting $p_i = P\{N = i\}$, $i = 1,2,3,\ldots$ we have from (1.1);

$$\phi(s) = \phi(cs)\ \{p_1 + p_2\ \phi(cs) + \ldots\}\ \text{for all}\ s>0\ \text{and some}\ c>0.$$

Since $0< \phi(s) <1$ and $\phi(s)$ is a decreasing function of $s$ one must have $\phi(cs) > \phi(s)$ and hence $cs<s$. This implies $0<c<1$. This appears to be true intuitively also as otherwise the distribution of the sum can grow out of proportion if the contribution of the individual components is not restricted. We will refer to $c$ as the scale parameter in the N-sum.

In the context of branching processes, when the PGF of $N$ is

$$P(t) = t/\{a - (a-1)t^k\}^{1/k},\ k>0\ \text{integer},\ a>1, \qquad (2.1)$$

Harris (1948) showed that the N-sum of gamma$(1,1/k)$ law with LT $1/[1+s]^{1/k}$ is of the type as gamma$(1,1/k)$. In fact the stability holds here for every $a>1$ as shown by Sandhya (1991a). When $k=1$ (2.1) becomes the PGF of the geometric law on $\{1, 2, \ldots\}$.

**Definition.2.1** The PGF (2.1) defines Harris$(a,k)$ law.

Now, let $\phi(u)$ be the CF that is N-sum stable w.r.t Harris$(a,k)$ for some $a>1$. Setting

$$\psi(u) = [1-\phi^k(u)] / \phi^k(u),\ \ \text{we have}\ \ \phi(s) = [1+\psi(u)]^{-1/k}.$$

Hence by the assumption that $\phi(u)$ is N-sum stable w.r.t Harris$(a,k)$ we have,

$$[1+\psi(u)]^{-1/k} = \frac{[1+\psi(cu)]^{-1/k}}{\{a-(a-1)/[1+\psi(cu)]\}^{1/k}} = [1+a\psi(cu)]^{-1/k}.$$

Hence $\psi(u) = a\psi(cu)$ for all $u \in \mathbf{R}$ and some $0<c<1<a$, that is $\psi(u)$ satisfies (1.2) with $b = c$. As the converse is easy we have proved:

**Theorem.2.1** A CF $\phi(u)$ is Harris$(a,k)$-sum stable for some $a>1$, if and only if $\phi(u) = 1/[1+\psi(u)]^{1/k}$, where $\psi(u)$ satisfies (1.2) with $\alpha \in (0,2]$, $b = c = a^{-1/\alpha}$, '$a$' and '$k$' being that in Harris$(a,k)$.

This theorem generalizes the result (Sandhya (1991a) and Theorem.3 of Lin (1994)) that a distribution on $\mathbf{R}$ is a geometric$(p)$-sum of its own type if and only if it is semi-$\alpha$-Laplace, which is a reformulation of Theorem.1 of Pillai (1985).

**Remark.2.2** The CF in Theorem.2.1 will be called as that of a generalized semi-$\alpha$-Laplace$(1/k)$ and the scale parameter $c$ in (1.1) is fixed by the parameter $b$ (or $a$) in this CF.

**Theorem.2.2** If $\phi(u)$ is N-sum stable w.r.t Harris$(a_1,k)$ and Harris$(a_2,k)$ such that $\ln a_1/ \ln a_2$ is irrational, then $\phi(u) = 1/[1+\lambda|u|^\alpha]^{1/k}$ (a CF mentioned in Devroye (1990)).

**Proof.** Follows from the property of $\psi(u)$ already mentioned in the introduction.

**Remark.2.3** In Theorem.2.1 suppose we assume that $\phi(u)$ is Harris$(a,k)$-sum stable for every $a>1$. Now by invoking Theorem.2.2 one is led to think that $\phi(u) = 1/[1+\lambda|u|^\alpha]^{1/k}$, because for every $a>1$ implies there exists $a_1>1$ and $a_2>1$ such that $\ln a_1/ \ln a_2$ is irrational. But this is not true in general.

In the context of geometric($p$)-sum stability it has been clarified by Lin (1998, Theorem.2) that $\phi(u)$ is geometric($p$)-sum stable for every $p \in (0,1)$ if and only if

$$\phi(u) = 1/[1+\lambda|u|^\alpha \, e^{-i\theta\,\mathrm{sgn}(u)}], \ |\theta| \leq \min(\pi\alpha/2, \pi-\pi\alpha/2).$$

(Of course the line of argument is different). This class of distributions called strictly geometric stable laws, is contained in semi-$\alpha$-Laplace laws (see the remark following Theorem.4.1 in Kozubowski and Rachev (1999)) and obviously contains the Linnik laws with CF $\phi(u) = 1/[1+\lambda|u|^\alpha]$ which is symmetric ($\theta = 0$). Since Lin's clarification is significant in the Harris($a,k$)-sum scheme also we add:

**Remark.2.4** In more generality one can show using equation (1.1) that the CF

$$\phi(u) = 1/[1+\lambda|u|^\alpha \, e^{-i\theta\,\mathrm{sgn}(u)}]^{1/k}, \ |\theta| \leq \min(\pi\alpha/2, \pi-\pi\alpha/2)$$

is Harris($a,k$)-sum stable for every $a>1$, thus establishing its N-sum stability property. Importantly, here the scale parameter $c$ in (1.1) is not fixed by any of the parameters in the CF. This class of distributions with the above CF is a subclass of the generalized Linnik($\alpha,\theta,\nu$) (GL) laws of Erdogan and Ostrovskii (1998), when $\nu = 1/k$.

From the PGF in (2.1) it follows that $\{a - (a-1)t^k\}^{-1/k}$, $k>0$ integer, $a>1$, is also a PGF. Hence by the N-sum stability property of GL($\alpha,\theta,1/k$) laws with CF $\phi(u)$, mentioned in Remark.2.4 it follows that

$$\phi(u) = \phi(cu)\{a - (a-1)[\phi(u)]^k\}^{-1/k} \text{ for every } c = a^{-1/\alpha} \text{ in } (0,1).$$

Thus we have, (see also Theorem.1 in Sandhya and Satheesh (1996));

**Theorem.2.3** GL($\alpha,\theta,1/k$) laws belong to class-L.

**Theorem.2.4** GL($\alpha,\theta,1/k$) laws are unimodal.

**Theorem.2.5** GL($\alpha,\theta,1/k$) laws are absolutely continuous.

From Theorem.2.1 and invoking the characterization of Linnik laws among semi-$\alpha$-Laplace laws in Sandhya and Satheesh (1996) by the class-L property, we have;

**Theorem.2.6** Generalized semi-$\alpha$-Laplace($1/k$) laws are in class-L if and only if it is that of a GL($\alpha,\theta,1/k$) law.

**Remark.2.5** It is known that semi stable laws are infinitely divisible (ID) and further (Ramachandran and Lau (1991), p.61) they are absolutely continuous. Now taking gamma($1,\beta$) mixtures of semi stable laws we have the CF $\phi(u) = 1/[1+\psi(u)]^\beta$, $\beta >0$ which we will refer to as that of generalized semi-$\alpha$-Laplace($\beta$) family. Being a continuous mixture they are also absolutely continuous. Though not a member of class-L in general, these distributions are ID as the mixing law here is ID.

The fact that generalized semi-$\alpha$-Laplace($\beta$) and GL($\alpha,\theta,\beta$) are mixtures of gamma($1,\beta$) laws readily suggest the corresponding stochastic representation which in turn provide a method to generate these variables from semi stable and strictly stable laws. Devroye (1990) has discussed this aspect for the CF $\phi(u) = 1/[1+\lambda|u|^\alpha]^\beta$. Erdogan and Ostrovskii (1998) have reached the conclusions in Theorems 2.4 and 2.5 (of course for the general $\nu >0$ and not just when $\nu = 1/k$) by different methods. Restricting the support of $X$ to $[0,\infty)$ we have the following corollary in terms of LTs.

**Corollary.2.1** A LT $\phi(s)$ is N-sum stable w.r.t. Harris$(a,k)$ law for some $a>1$, if and only if $\phi(s) = [1+\psi(s)]^{-1/k}$, where $\psi(s)$ satisfies (1.2) with $\alpha \in (0,1]$, $b = c$, '$a$' and '$k$' being that in Harris$(a,k)$. Further, the LT $1/(1+s^\alpha)^{1/k}$ is Harris$(a,k)$-sum stable for every $a > 1$.

This is a generalization of the semi Mittag-Leffler (SML$(a,b)$) law given by the LT $1/[1+\psi(s)]$ where $\psi(s)$ satisfies (1.2) which was introduced and characterized by Sandhya (1991b) as the family of distributions invariant under p-thinning for some $p$, thus generalizing the result of Renyi. ML laws with LT $1/(1+s^\alpha)$, $0<\alpha \leq 1$, are known to be geometrically stable for every $p$ in $(0,1)$. Since SML$(a,b)$ laws are ID we can have a generalized SML$(a,b)$ laws with LT $\phi(s) = [1+\psi(s)]^{-\beta}$, $\beta > 0$. But rather curiously this in general is not N-sum stable w.r.t Harris law, which will be proved in Theorem.4.1.

## 3. Stability of integer valued random variables

The following lemma is instrumental in developing (i) PGFs representing the discrete analogue of the LTs of continuous laws on $[0,\infty)$, (ii) the discrete analogue of the notion of distributions of the same type and (iii) stability of integer valued r.vs in the N-sum scheme.

**Lemma.3.1** If $\phi(s)$ is a LT, then $Q(s) = \phi(1-s)$, $0<s<1$ is a PGF.

**Proof.** As $\phi(s)$ is a LT, it is completely monotone and $\phi(o) = 1$. Hence $Q(s)$ is absolutely monotone and $Q(1) = 1$. Hence by Feller (1971, 223) $Q(s)$ is a PGF.

From this Lemma we have the following PGFs as discrete analogue of their continuous counterparts.

(i) The PGF of discrete stable laws of Steutel and van Harn (1979) is,

$$Q(s) = \exp\{-\lambda(1-s)^\alpha\}, \lambda>0, 0<\alpha \leq 1; \text{ when } \phi(s) = \exp\{-\lambda s^\alpha\}. \tag{3.1}$$

(ii) The PGF of discrete ML laws of Pillai and Jayakumar (1995) is,

$$Q(s) = 1/\{1+\lambda(1-s)^\alpha\}, \lambda>0, 0<\alpha \leq 1; \text{ when } \phi(s) = 1/\{1+\lambda s^\alpha\}. \tag{3.2}$$

(iii) The PGF of discrete Linnik laws of Pakes (c.f Devroye(1993)) is,

$$Q(s) = 1/\{1+\lambda(1-s)^\alpha\}^\beta, 0<\alpha \leq 1, \lambda>0, \beta>0; \text{ when } \phi(s) = 1/\{1+\lambda s^\alpha\}^\beta. \tag{3.3}$$

(iv) The PGF $Q(s) = 1/[1+\psi(1-s)]^\beta$ where $\psi$ satisfies (2) with $\alpha \in (0,1]$ and $\beta > 0$, (3.4)

describes the discrete generalized SML$(a,b)$ laws; when $\phi(s) = [1+\psi(s)]^{-\beta}$, $\beta > 0$.

It is worth mentioning here that the geneses of the first three discrete laws as given by these authors were not as simple and straightforward as conceived here.

Two LTs $\phi_1$ and $\phi_2$ are of the same type if $\phi_1(s) = \phi_2(cs)$ for all $s>0$ and some $c>0$. Setting $Q_1(s) = \phi_1(1-s)$ and $Q_2(s) = \phi_2(1-s)$ (by Lemma.3.1) and $c<1$, this is reflected in $Q_1$ and $Q_2$ as $Q_1(1-s) = Q_2(1-cs)$. Equivalently, putting $t = 1-s$, $Q_1(t) = Q_2(1-c+ct)$, for all $0<t<1$. Notice that in the N-sum scheme $c<1$ by Remark.2.1. These justify the notion of D-type, (distributions of the same type for discrete laws) below.

**Definition.3.1** Two PGFs $Q_1(s)$ and $Q_2(s)$ are of the same D-type if $Q_1(1-s) = Q_2(1-cs)$, for all $0<s<1$, or equivalently, $Q_1(u) = Q_2(1-c+cu)$ for all $0<u<1$, and some $0<c<1$.

Using the second representation of PGFs of the same D-type the definition of discrete class-L laws of Steutel and van Harn (1979) follows directly and consequently the PGF constructed (by Lemma.3.1) from a LT in class-L is in discrete class-L. Thus:

**Theorem.3.1** The discrete Linnik laws given by (3.3) with $\beta = 1/k$, where $k$ is a positive integer, is in discrete class-L.

**Theorem.3.2** The PGF in (3.4) is in discrete class-L if and only if it is a discrete Linnik law given by (3.3) with $\beta = 1/k$, where $k>0$ is an integer (the discrete analogue of Theorem.2.6).

From Definition.3.1, the stochastic representation corresponding to two integer valued r.vs $X$ and $Y$ of the same D-type is $X \stackrel{d}{=} \sum_{i=1}^{Y} Z_i$, where $\{Z_i\}$ are i.i.d Bernoulli($c$) variables independent of $Y$ with $P\{Z_i = 0\} = 1-c$. Thus using $c \circ X = \sum_{i=1}^{X} Z_i$, as the discrete analogue of $cX$ done in Steutel and van Harn (1979) is also justified.

The discrete stable laws also are stable under ordinary summation quite like the stable (continuous) laws (Feller (1971, p.448)). This follows from;

$$[\exp\{-\lambda[k^{-1/\alpha}(1-s)]^{\alpha}\}]^k = [\exp\{-\lambda k^{-1}(1-s)^{\alpha}\}]^k = \exp\{-\lambda(1-s)^{\alpha}\}$$

Here $c = k^{-1/\alpha}$ is the scale parameter. Steutel and van Harn (1979) define discrete stable laws by PGFs satisfying,

$$Q(s) = Q[1-c(1-s)] \cdot Q[1-(1-c^{\alpha})^{1/\alpha}(1-s)], \quad 0<c<1 \tag{3.5}$$

and arrived at the form in (3.1). Here 'discrete stability' is described in the sense;

$$X \stackrel{d}{=} c \circ X_1 + (1-c^{\alpha})^{1/\alpha} \circ X_2.$$

Rao and Shanbhag (1994, p.160) arrives at (3.1) from (3.5) using a different approach. Here we have obtained (3.1) directly from the LT of the continuous case and further shown that the PGF thus obtained is 'stable' under ordinary summation.

In the discrete set up corresponding to equation (1.1), the description of stability of N-sums can now be described by; with $Q(s)$ the PGF of $X$;

$$P[Q(1-c+cs)] = Q(s) \text{ for all } s \in (0,1) \text{ and some } 0<c<1. \tag{3.6}$$

We thus have the following results in the discrete domain invoking Lemma.3.1 and (3.6).

**Theorem.3.3** A PGF $Q(s)$ is Harris($a,k$)-sum stable for some $a>1$, if and only if $Q(s) = 1/[1+\psi(1-s)]^{1/k}$, where $\psi(1-s)$ satisfies (1.2), $b = c = a^{-1/\alpha}$, '$a$' and '$k$' being that in Harris($a,k$).

**Corollary.3.1** A PGF $Q(s)$ is a geometric($p$)-sum of its own D-type if and only if it is discrete SML($a,b$) with PGF $Q(s) = 1/[1+\psi(1-s)]$, $a = 1/p$ and $c = b$.

**Corollary.3.2** A PGF $Q(s)$ is a geometric($p$)-sum of its own D-type for every $p$ if and only if it is discrete ML with PGF $Q(s) = 1/[1+\lambda(1-s)^{\alpha}]$, $\lambda>0$.

This property of discrete ML laws has not been mentioned in Pillai and Jayakumar (1995). Also, if we consider in the context of p-thinning, a renewal counting process these corollaries mean that, a renewal counting process is invariant under p-thinning for some $p$ if and only if it is discrete SML($a,b$), and for every $p$ if and only if it is discrete ML which are discrete analogues of properties of SML laws in Sandhya(1991b).

Notice that when a continuous r.v $X$ is N-sum stable, $P\{N=0\}>0$ is clearly impossible because in such a case the N-sum $S_N$ can vanish with $P\{N=0\}$ while $X$ cannot. But, considering the PGFs $Q(s) = 1- \delta(1-s)^v$, $0<\delta<1$, $0<v<1$ and $P(s) = 1-\lambda(1-s)$, $0<\lambda<1$, $P[Q(s)] = 1- \lambda\delta(1-s)^v$, is again of the same D-type. Thus in the discrete set up $P\{N=0\}>0$ is possible. Here $Q(s)$ is a Sibuya($v$)-sum (with PGF $1-(1-s)^v$) of Bernoulli($\delta^{1/v}$) laws.

## 4. Identifying the compounding r.v $N$ for N-sum stability

From Sandhya (1991a, 1996) we have the following result on N-infinitely divisible laws: A necessary and sufficient condition that a r.v $X$ with LT $\phi$ is N-ID is that there exists a LT $\phi_m$ such that $\phi[\phi_m^{-1}]$ is a PGF for every $m$ in the parameter space of $m$. Sandhya (1996) has used a similar relation to identify $N$ for the N-sum stability of Mittag-Leffler laws. We now generalize this. In the N-sum stability scheme, when $X$ is non-negative, the key equation is;

$$Q[\phi(cs)] = \phi(s), \text{ for all } s>0 \text{ and some } c\in (0,1). \tag{4.1}$$

Setting $\phi_c(s) = \phi(cs)$ we have from (4.1) that

$$Q[\phi_c(s)] = \phi(s) \text{ for all } s>0 \text{ and some } c\in (0,1).$$

Since LTs are decreasing in $(0,\infty)$ the inverse function $\phi^{-1}$ exists uniquely. Hence setting

$$s = \phi_c^{-1}(t) \text{ for } 0<t<1, \text{ we have;}$$

$$Q(t) = \phi[\phi_c^{-1}(t)] \text{ for some } c\in (0,1). \tag{4.2}$$

Thus the PGF $Q$ of the compounding law is identified uniquely from the LT $\phi$ under the stability condition (4.1). Similarly, since PGFs are increasing in (0,1), their inverses exist uniquely, and hence from (3.6), when $s$ equals $Q_c^{-1}(t)$, $0<t<1$ where $Q_c(s) = Q(1-c+cs)$ we have;

$$P(t) = Q[Q_c^{-1}(t)], \text{ for some } 0<c<1. \tag{4.3}$$

One may as well conceive (4.2) and (4.3) as definitions of N-sum stability. (After the completion of this work we have come across Bunge (1996) who has considered a relation similar to (4.2), not as a method to identify $N$ for the N-sum stability of $X$, but for generating PGFs from LTs). Also, (4.2) and (4.3) suggest that $c$ appears as a parameter in the PGF. The range of the values of $c$ and in turn its relation to other parameters in $P, \phi$ or $Q$ in the above scheme remains to be determined. For this purpose, we establish the following Lemma. Incidentally it also replaces certain topological semi-group arguments in Bunge (1996) to arrive at the PGF.

**Lemma.4.1** If $P(s)$ is a PGF, then $P(s^u)$ is a PGF if and only if $u>0$ is an integer.

**Proof.** $P(s^u)$ is a PGF if and only if it can be represented as a power series in $s$ (Feller, (1971), p.223). When $u>0$ is not an integer this is not possible, proving the assertion.

We begin with a general result obtained by using (4.2) and Lemma.4.1.

**Theorem.4.1** The generalized SML($a,b,\beta$) law with LT $\phi(s) = [1+\psi(s)]^{-\beta}$, $\beta>0$ is N-sum stable if and only if $N$ is Harris($a,k$), $\beta = 1/k$, $c = b$, '$a$' and '$k$' being same in both.

**Proof.** We have $\phi(s) = [1+\psi(cs)]^{-\beta} = [1+a\psi(bs)]^{-\beta}$. Now, $\phi_c(s) = [1+\psi(cs)]^{-\beta}$ and

$$\phi_c^{-1}(t) = \frac{1}{c}\psi^{-1}\left(\frac{1-t^{1/\beta}}{t^{1/\beta}}\right). \quad \text{Hence}$$

$$\phi[\phi_c^{-1}(t)] = \left\{1+a\psi\left[b\frac{1}{c}\psi^{-1}(v)\right]\right\}^{-1}, \quad v = \frac{1-t^{1/\beta}}{t^{1/\beta}}$$

$$= 1/(1+av), \text{ when } c = b$$

$$= \left\{\frac{\omega t^{1/\beta}}{1-(1-\omega)t^{1/\beta}}\right\}^{\beta}, \quad \omega = 1/a.$$

Since $t\{\omega/[1-(1-\omega)t]\}^{\beta}$ is a PGF for $0<\omega<1$ and $\beta>0$, the above function is a PGF only when further $1/\beta = k$ a positive integer, by invoking Lemma.4.1. Hence $N$ is Harris$(a,k)$ where $a = 1/\omega$. The converse is clear and the proof is complete.

**Corollary.4.1** Setting $\psi(s) = s^{\alpha}$ we have a positive generalized Linnik$(\alpha,0,\beta)$ law (generalized ML). This is N-sum stable if and only if $N$ is Harris$(a,k)$, $\beta k = 1$ and $ac^{\alpha} = 1$.

**Corollary.4.2** Setting $k = 1$ we have; SML$(a,b)$ law is N-sum stable if and only if $N$ is geometric$(p)$, $c = b$ and $p = 1/a = c^{\alpha}$. (This is an example in Bunge (1996)).

**Corollary.4.3** Setting $\psi(s) = s^{\alpha}$, and $k = 1$ we have; ML law is N-sum stable if and only if $N$ is geometric with parameter $p = c^{\alpha}$, $0<c<1$. (The result of Sandhya (1996)).

**Corollary.4.4** Further putting $\alpha = 1$, ML law reduces to the exponential law and we have; the exponential law is N-sum stable if and only if $N$ is geometric with parameter $p = c$. (The result of Cinlar and Agnew (1968)).

**Corollary.4.5** Setting $\psi(s) = s$ we have; gamma$(1,\beta)$ law is N-sum stable if and only if $N$ is Harris$(a,k)$, $\beta = 1/k$ and $c = 1/a$. (Another example in Bunge (1996)).

**Remark.4.1** A point to be stressed here is that the generalized SML$(a,b,\beta)$ law, and in particular the gamma$(1,\beta)$ law, cannot be N-sum stable unless $\beta = 1/k$, where $k$ is a positive integer. (See the last line of Section.2 in this regard).

Considering the LT of a semi stable (SS$(a,b)$) law we have $\phi(s) = \exp\{-\psi(s)\} = \exp\{-a\psi(bs)\}$. With $\phi_c(s) = \exp\{-\psi(cs)\}$ we have $\phi_c^{-1}(t) = \frac{1}{c}\psi^{-1}\left(\ln\frac{1}{t}\right)$. When $c = b$, $\phi[\phi_c^{-1}(t)] = t^a$ which is a PGF only when $a>1$ is an integer say $k$. Since $kc^{\alpha} = 1$ (as $ab^{\alpha} = 1$), $0<\alpha\leq 1$ we have proved:

**Theorem.4.2** A SS$(a,b)$ law is N-sum stable if and only if $a>1$ is an integer, say $k$, $N$ is degenerate at $k$ and $c = b = k^{-1/\alpha}$.

When $\psi(s) = s^{\alpha}$, SS$(a,b)$ law becomes a stable$(\alpha)$ law (Feller, 1971, 448) and following a similar line of argument we have; $Q[Q_c^{-1}(t)] = t^{1/c^{\alpha}}$ which is a PGF only when $1/c^{\alpha} = k>0$ integer, and we have proved,

**Theorem.4.3** A stable$(\alpha)$ law is N-sum stable if and only if $P\{N=k\} = 1$ for any (arbitrary) $k\geq 1$ integer and $c = k^{-1/\alpha}$.

Feller (1971) had shown that if $X$ is stable($\alpha$), then $X \stackrel{d}{=} n^{-1/\alpha}[X_1+\ldots+X_n]$ for every $n$. Here we have shown that when $X$ is stable($\alpha$) a fixed sum alone imparts stability in the $N$-sum scheme and $c = n^{-1/\alpha}$ for every $n$, the number of components in the sum. Notice that for a SS($k,b$) law $k>1$ integer fixes the number of components in the sum and it cannot vary, where as in the case of stable($\alpha$) laws the number of components $n$ may vary, since a corresponding change in the scale parameter $c = n^{-1/\alpha}$ will guarantee stability. Hence if a semi stable law shows stability w.r.t every $k$ then it must be stable($\alpha$). This conclusion is supported also by the fact that if this is true for every positive integer $n$, then there exists two integers say, $n_1$ and $n_2$ such that $\ln n_1 / \ln n_2$ is irrational and consequently $\psi(s) = \lambda s^\alpha$ for a positive constant $\lambda$ (see the description of equation (1.2)). Intuitively it thus appears that the condition

$$X \stackrel{d}{=} N^{-1/\alpha}\{X_1 + \ldots + X_N\}$$

where $N$ is a positive integer valued r.v, characterizes the stable($\alpha$) law. This has been proved in Ramachandran and Lau (1991, p.88).

**Remark.4.2** The method presented here can be used to identify $N$ even when $X$ is $\nu$-stable (and not necessarily positive) as conceived in Gnedenko and Korolev (1996), Bunge (1996) and Kozubowski and Panorska (1996). This is because these CFs are of the form $\phi(\psi(u))$ where $\phi$ is a LT and $\psi(u) = \lambda|u|^\alpha e^{-i\theta\mathrm{sgn}(u)}$, $|\theta| \leq \min(\pi\alpha/2, \pi-\pi\alpha/2)$ and using (4.2) we have $\phi\{\phi_c^{-1}[\phi_c(\psi(u))]\} = \phi\{\psi(u)\}$. Also when $\psi(u)$ satisfies (1.2) we have the CFs of the corresponding $\nu$-semi-stable laws and under the condition $c = b$ the method still works.

Invoking Lemma.3.1 we have discrete analogues of gamma, ML their generalizations and results analogous to those obtained in this section with proofs on similar lines. Considering the importance of the discrete setup we only state the following general results, omitting the proofs.

**Theorem.4.4** The PGF $Q(s) = [1+\psi(1-s)]^{-\beta}$, where $\psi(1-s)$ satisfies (1.2) is N-sum stable if and only if $N$ is Harris($a,k$), $\beta = 1/k$, $k>0$ integer, $c = b$, '$a$' and '$k$' being same in both.

**Theorem.4.5** The discrete SS($a,b$) law with PGF $\exp\{-\psi(1-s)\}$ is N-sum stable if and only if the parameter $a>1$ is an integer say $k$, $N$ is degenerate at $k$ and $c = b = k^{-1/\alpha}$.

Also, the geometric law on $\{0,1,2,\ldots\}$ and discrete analogues of SML and ML laws are stable only w.r.t a geometric-sum, and the discrete analogue of gamma$(1,\beta)$ is stable only w.r.t a Harris-sum and further $\beta = 1/k$ for a positive integer $k$.

**Acknowledgement.** The authors wish to thank the referee for bringing to their notice the reference Lin (1998), its relevance and suggesting to include Remark.2.3.